\input amstex
\documentstyle{amsppt}
\magnification=\magstep1 \baselineskip=12pt \hsize=6truein
\vsize=8truein

\topmatter

\title
Compactness for conformal metrics with Constant $Q$ curvature on
locally conformally flat manifolds
\endtitle

\author Jie Qing and David Raske
\endauthor

\leftheadtext{Compactness} \rightheadtext{Jie Qing and David
Raske}

\address Jie Qing, Department of Mathematics, UC Santa Cruz, Santa Cruz, CA
95064 \endaddress \email qing{\@}ucsc.edu\endemail
\address David Raske, Department of Mathematics, UC Santa Cruz, Santa Cruz, CA
95064 \endaddress \email gadfly{\@}math.ucsc.edu\endemail

\abstract In this note we study the conformal metrics of constant
$Q$ curvature on closed locally conformally flat manifolds. We
prove that for a closed locally conformally flat manifold of
dimension $n\geq 5$ and with Poincar\"{e} exponent less than
$\frac {n-4}2$, the set of conformal metrics of positive constant
$Q$ and positive scalar curvature is compact in the $C^\infty$
topology.
\endabstract
\endtopmatter

\document
\noindent{\bf 1. Introduction}\vskip 0.1in

Let $(M, g)$ be a compact manifold without boundary of dimension
higher than $4$. Let
$$
Q[g] = - \frac {n-4}{4(n-1)} \Delta R + \frac {(n-4)(n^3 - 4n^2 +
16n - 16)}{16(n-1)^2(n-2)^2}R^2 - \frac {2(n-4)}{(n-2)^2}|Ric|^2
\tag 1.1
$$
be the so-called $Q$-curvature, where $R$ is the scalar curvature,
$Ric$ is the Ricci curvature. And let
$$
P[g] = (-\Delta)^2 - \text{div}_g((\frac {(n-2)^2 +
4}{2(n-1)(n-2)} Rg - \frac 4{n-2}Ric_g)d) + Q[g] \tag 1.2
$$
be the so-called Paneitz-Branson operator. It is known that
$$
P[g] w = Q[g_w] w^{\frac {n+4}{n-4}} \tag 1.3
$$
which is called the Paneitz-Branson equation, where $g_w =
w^{\frac 4{n-4}}g$ ({\it cf}. [P] [Br] [DHL] [DMA] ). We consider
the equation (1.3) as a fourth order analogue of the well-known
scalar curvature equation
$$
L[g] u = R[g_u] u^{\frac {n+2}{n-2}}, \tag 1.4
$$
where
$$
L[g] = - \frac {4(n-1)}{n-2}\Delta + R \tag 1.5
$$
is the so-called conformal Laplacian and $g_u = u^{\frac
4{n-2}}g$. The well-known Yamabe problem in conformal geometry is
to find a metric, in a given class of conformal metrics, which is
of constant scalar curvature, {\it i.e.} to solve
$$
L[g] u =Y u^{\frac {n+2}{n-2}}
$$
on a given manifold $(M, g)$ for some positive function $u$ and a
constant $Y$. The affirmative resolution to the Yamabe problem was
given in [Sc1] after other notable works [Ya] [Tr] [Au].

Recently, there are more and more interests in using higher order
partial differential equations in the study of conformal geometry,
particularly after some successes in dimension 4 in [CY]
(references therein). One major hurdle that stands in the way to
this higher order approach is the lack of the maximum principles
for the higher order partial differential equations like (1.3).
One, for example, would like to say that a nonnegative solution
$w$ to (1.3) is necessarily positive, and therefore gives rise to
a metric $g_w = w^{\frac 4{n-4}}g$.

In this note we consider closed locally conformally flat manifolds
with positive Yamabe constant of dimension higher than $4$. By
work in [SY], we know that such manifolds are all Kleinian in the
sense that, there is a Kleinian group $\Gamma$ such that $(M,
[g])$ is conformally equivalent to $\Omega(\Gamma)/\Gamma$, where
$\Omega(\Gamma)$ is set of ordinary points of the Kleinian group
$\Gamma$ on $S^n$. The developing map is a conformal
diffeomorphism from the universal cover $\tilde{M}$ of $M$ to
$\Omega(\Gamma) \subset S^n$. The Kleinian group $\Gamma$ is the
so-called holonomy representation of the fundamental group of $M$.
Hence we may turn to study solutions to
$$
(-\Delta)^2 \hat w = Q(\hat g_w) \hat w^\frac {n+4}{n-4} \tag 1.6
$$
on a domain $\hat \Omega \subset R^n$ instead of solutions to the
Paneitz-Branson equation (1.3) on $M$. Given a Kleinian group
$\Gamma$, the so-called Poincar\'{e} exponent is defined as:
$$
\delta(\Gamma) = \inf\{\delta > 0: \sum_{\gamma\in \Gamma}
|\gamma'(x)|^\delta < \infty, \forall x\in S^n\}.
$$
Inspired by the work of Wei and Xu in [WX], where they managed to
use some maximum principle and the moving plane method to classify
all positive solutions to (1.6) on entire $R^n$ with $Q(\hat g_w)$
a positive constant, we have

\proclaim{Theorem 1.1} Suppose that $(M^n, g)$ is a closed locally
conformally flat manifold with positive Yamabe constant and its
Poincar\'{e} exponent less than $\frac {n-4}2$. Then any
nontrivial nonnegative solution to
$$
P[g] u = Q u^\frac {n+4}{n-4}
$$
on $M$ with positive constant $Q$ has to be strictly positive on
$M$.
\endproclaim

The assumption on Poincar\'{e} exponent gives us the integrability
condition to use the idea from [WX]. We also adopt the nice and
geometric use of the moving plane method of Schoen from [Sc2] and
prove a local convexity result analogous to the one in [Sc2]. One
key point is to use the Hopf maximum principle for elliptic
systems instead (cf. [Q]).

\proclaim{Lemma 1.2} Let $(\Omega, \tilde g)$ be the universal
cover of a closed manifold $(M^n, g)$ where $n> 4$, $\Omega$ is
not $S^n$, $\tilde g = v^\frac 4{n-4} g_1$ and $g_1$ is the
standard round metric on $S^n$. And suppose that $Q[\tilde g]$ is
a positive constant and that the scalar curvature $R[\tilde g]$ is
positive. Then any small round ball $B\subset \Omega \subset S^n$
is geodesically convex with respect to the metric $\tilde g$.
\endproclaim

Combining Lemma 1.2 in the above with the classification given in
[WX], we therefore derive our $L^\infty$ estimate:

\proclaim{Theorem 1.3} Let $(M, g)$ be a closed locally
conformally flat manifold of dimension greater than $4$ with
positive Yamabe constant. And suppose that $(M, g)$ is not
conformal to $(S^n, g_1)$. Then there exists $C >0$ such that, if
$g_u = u^\frac 4{n-4}g$ is of positive scalar curvature and
constant $Q$-curvature $1$, then
$$
\|u\|_{L^\infty(M)} \leq C. \tag 1.7
$$
\endproclaim

Similar to the study of Yamabe problem, we consider
$$
P(M, [g]) = \inf_{u^\frac 4{n-4}g\in [g]} \frac {\int_M
(Qdv)[u^\frac 4{n-4}g]}{(\int_M dv[u^\frac 4{n-4}g])^\frac
{n-4}n}. \tag 1.8
$$
We will refer to $P(M, [g])$ as the Paneitz constant in this paper. Clearly the
Paneitz constant $P(M, [g])$ is a conformal invariant. As a
consequence of Theorem 1.1 and Theorem 1.3 in the above, here is
our main theorem:

\proclaim{Theorem 1.4} Let $(M^n, g)$ be a closed, locally
conformally flat manifold of dimension $n$ greater than $4$ with
positive Yamabe constant and Poincar\'{e} exponent less than
$\frac {n-4}2$. Suppose that $(M, g)$ is not conformal to $(S^n,
g_1)$ and that the Paneitz constant $P(M, [g])$ is positive. Then
the set of metrics conformal to $g$ on $M$ with positive scalar
curvature and constant $Q$-curvature $1$ is compact in $C^\infty$
topology. That is, there exists $C_k >0$ for each $k\in N$ such
that, if $g_u = u^\frac 4{n-4}g$ is of positive scalar curvature
and constant $Q$-curvature $1$, then
$$
\|u\|_{C^k(M)} + \|\frac 1u \|_{C^k(M)} \leq C_k. \tag 1.9
$$
\endproclaim

We would like to mention that it was proven in [CHY] [G] that the
Hausdorff dimension of the limit set of the Kleinian group
associated with a closed, locally conformally flat manifold $(M, g)$
of dimension higher than $4$ with positive scalar curvature and
positive $\sigma_2$-curvature is less than $\frac {n-4}2$, which in
turn implies the Poincar\'{e} exponent is less than $\frac {n-4}2$.
$\sigma_2$-curvature is the second symmetric function of the
eigenvalues of Weyl-Schouten tensor $R_{ij} - \frac
1{2(n-1)}Rg_{ij}$. It is certainly conceivable that the positivity
of both the Paneitz constant and Yamabe constant should imply the
Poincar\'{e} exponent is less than $\frac {n-4}2$.

Finally we like to mention that after submission of this paper we
learnt the announcement of Hebey and Robert [HR], where they had
obtained some results similar to the ones in this paper. One
difference is that our compactness is for the metrics, while their
compactness is for the nonnegative solutions. In other words, we
have Theorem 1.1 and the positive lower bound in this paper. Our
approach is to use the moving plane method to obtain the geodesic
convexity of round balls and exclude the possibility of blow-ups,
while the analysis in [HR] is based on the study of the asymptotic
analytic behavior of solutions near the blow-up points and uses the
strong positivity assumption on the Green function.

The organization of the note is as follows: in Section 2 we will
recall more about Kleinian groups and works in [SY] to prove Theorem
1.1. In Section 3 we will adopt the moving plane approach from [Sc2]
to prove Lemma 1.2. Finally in Section 4 we use the blow-up method
to obtain $L^\infty$ estimates and prove our main theorem.

\vskip 0.1in\noindent {\bf Acknowledgement} \quad We would like to
thank Professor Alice Chang and Professor Paul Yang for their
interests in this work.

\vskip0.1in \noindent{\bf 2. Positivity of solutions to
Paneitz-Branson equations}\vskip 0.1in

Suppose that $(M^n, g)$ is a closed locally conformally flat
manifold with positive Yamabe constant. Then by works in [SY] we
know that the developing map from the universal cover $(\tilde M,
\tilde g)$ into $(S^n, g_1)$ is a conformal embedding, where $g_1$
is the standard round metric on $S^n$. And the deck transformation
group of the universal covering becomes a discrete subgroup of the
group of conformal transformation of $S^n$, which is called a
Kleinian group. This Kleinian group $\Gamma$ is said to be the
holonomy representation of the fundamental group of $M$. The image
of the developing map from the universal cover $\tilde M$ is the
set $\Omega(\Gamma)$ of ordinary points for the Kleinian group
$\Gamma$, and $(M, [g])$ is conformally equivalent to
$\Omega(\Gamma)/\Gamma$, which is called the Kleinian manifold
associated with the group $\Gamma$. Therefore we may consider
$(\Omega(\Gamma), \tilde g)$ as the universal cover for $(M, g)$
and $\tilde g = \tilde \eta^{\frac 4{n-4}}g_1$ . Hence
$$
P[g_1]\tilde \eta = Q[\tilde g]\tilde \eta^{\frac {n+4}{n-4}}.
\tag 2.1
$$
Using stereographic projection
$$
\psi: \hat\Omega \subset R^n \longrightarrow \Omega
(\Gamma)\subset S^n \tag 2.2
$$
with respect to some point in $\Omega(\Gamma)\subset S^n$, we may
write $\hat g  =  \psi^* \tilde g= \hat \eta^{\frac 4{n-4}}g_0$,
where $g_0$ is the Euclidean metric and $\hat \eta = (\tilde\eta
\circ \psi)(\frac 2{1+|x|^2})^\frac {n-4}2$, and
$$
P[g_0]\hat \eta = Q[\hat g]\hat \eta^{\frac {n+4}{n-4}}. \tag 2.3
$$
Notice that $P[g_0] = (-\Delta)^2$. Suppose that $g_u = u^\frac
4{n-4} g$ is another metric in $[g]$ on $M$ and $\phi:
\Omega(\Gamma) \rightarrow M$ is the covering map, which is the
composition of the developing map and the universal covering. Then
$$
\hat g_u  = \hat u^\frac 4{n-4} g_0, \tag 2.4
$$
where $\hat u = (u \circ \phi \circ \psi) \hat \eta$. Then
$$
P[g] u = Q[g_u]u^\frac {n+4}{n-4} \quad \text{on $M$} \tag 2.5
$$
and
$$
(-\Delta)^2 \hat u = Q[\hat g_u] \hat u^\frac {n+4}{n-4} \quad
\text{in $\hat \Omega\subset R^n$} \tag 2.6
$$
where $Q[\hat g_u] = Q[g_u]\circ \phi\circ \psi$. Due to the
invariant property of the above equations, one may easily see
that, for a nonnegative function $v$ satisfying
$$
P[g] v = Q v^\frac {n+4}{n-4} \quad \text{on $M$} \tag 2.7
$$
for some function $Q$ on $M$, $\hat v = (v\circ \phi\circ \psi)
\hat \eta$ satisfies
$$
(-\Delta)^2 \hat v = \hat Q \hat v^\frac {n+4}{n-4} \quad \text{in
$\hat \Omega\subset R^n$} \tag 2.8
$$
for the function $\hat Q = Q\circ \phi\circ \psi$ on $\hat
\Omega$.

One very important quantity of a Kleinian group $\Gamma$ is the
so-called Poincar\'{e} exponent
$$
\delta(\Gamma) = \inf\{\delta > 0: \sum_{\gamma\in \Gamma}
|\gamma'(x)|^\delta < \infty, \forall x\in S^n\}. \tag 2.9
$$
Another very important quantity is the Hausdorff dimension
$d(\Gamma)$ of the limit set $L(\Gamma)$. Due to a theorem of
Patterson and Sullivan one knows that $\delta(\Gamma) = d(\Gamma)$
when $\Gamma$ is geometrically finite. In our cases, the Kleinian
group associated with a closed, locally conformally flat manifold
of positive Yamabe constant is always geometrically finite (cf.
[CQY2]). Now we are ready to state and prove an interesting
property of a Kleinian manifold.

\proclaim{Lemma 2.1} Suppose that $(M, [g])$ is a Kleinian
manifold $\Omega(\Gamma)/\Gamma$ with $\delta(\Gamma) < \frac
{n-4}2$. And suppose that $v$ is a nonnegative function on $M$ and
$\tilde v = (v\circ \phi)\tilde \eta$. Then
$$
\int_{S^n} \tilde v^\frac {n+4}{n-4} dv_{g_1} < \infty. \tag 2.10
$$
\endproclaim
\demo{Proof} We need only to show that
$$
\int_{\Omega(\Gamma)}\tilde v^\frac {n+4}{n-4} dv_{g_1} < \infty.
\tag 2.11
$$
Because the Hausdorff dimension $d(\Gamma)$ of $S^n\setminus
\Omega(\Gamma)$ is less than $\frac {n-4}2$ by the theorem of
Patterson and Sullivan. Notice that, if we let $F$ be a
fundamental domain of $\Gamma$,
$$
\int_{\Omega(\Gamma)}\tilde v^\frac {n+4}{n-4} dv_{g_1} =
\sum_{\gamma\in \Gamma} \int_{\gamma(F)} \tilde v^\frac {n+4}{n-4}
dv_{g_1}. \tag 2.12
$$
Since $\gamma$ is an isometry for $g$, we have
$$
\tilde \eta^\frac 4{n-4}(x) g_1 = \gamma^* (\tilde \eta^\frac
4{n-4}(x) g_1) = \tilde \eta ^\frac 4{n-4}(\gamma(x)) |\gamma'|^2
g_1 \tag 2.14
$$
and
$$
\tilde v (\gamma(x)) = v \circ \phi(\gamma(x)) \tilde
\eta(\gamma(x)) = \tilde v (x) |\gamma'(x)|^{-\frac {n-4}2}. \tag
2.15
$$
Hence
$$
\int_{\gamma(F)}\tilde v^\frac {n+4}{n-4} dv_{g_1} = \int_{F}
\tilde v (x)^\frac {n+4}{n-4} |\gamma'(x)|^{-\frac {n+4}2}
|\gamma'(x)|^n dv_{g_1} = \int_{F} \tilde v (x)^\frac {n+4}{n-4}
|\gamma'(x)|^\frac {n-4}2 dv_{g_1}. \tag 2.16
$$
Therefore
$$
\int_{\Omega(\Gamma)}\tilde v^\frac {n+4}{n-4} dv_{g_1} = \int_F
(\sum_{\gamma\in \Gamma}|\gamma'(x)|^\frac {n-4}2) \tilde v
(x)^\frac {n+4}{n-4}dv_{g_1}, \tag 2.17
$$
which is finite when the Poincar\'{e} exponent, $\delta(\Gamma) <
\frac {n-4}2$. Thus the proof is completed.
\enddemo

\proclaim{Corollary 2.2} Suppose that $(M, [g])$ is a closed
Kleinian manifold $\Omega(\Gamma)/\Gamma$ with $\delta(\Gamma) <
\frac {n-4}2$. And suppose that $v$ is a nonnegative function on
$M$ and $\hat v = (v \circ \phi\circ \psi)\hat \eta$. Then
$$
\int_{U} \hat v^\frac {n+4}{n-4} dv_{g_0} < \infty \tag 2.18
$$
for any bounded domain $U\subset R^n$.
\endproclaim
\demo{Proof} Notice that
$$
\hat v = \tilde v \circ \psi (\frac 2{1+|x|^2})^\frac {n-4}2. \tag
2.19
$$
Hence
$$
\int_{U} \hat v^\frac {n+4}{n-4} dv_{g_0}  = \int_{\psi (U)}
\tilde v^\frac {n+4}{n-4} (\frac 2{1+|x|^2})^{- \frac {n-4}2}
dv_{g_1}. \tag 2.20
$$
Therefore the corollary is proved.
\enddemo

Next we state and prove a theorem on the positivity of solutions
to the Paneitz-Branson equation on $M$.

\proclaim{Theorem 2.3} Suppose that $(M, [g])$ is a closed
Kleinian manifold $\Omega(\Gamma)/\Gamma$ with $\delta(\Gamma) <
\frac {n-4}2$. And suppose that $v$ is a smooth nonnegative
function on $M$ satisfying
$$
P[g] v = Q v^\frac {n+4}{n-4} \quad \text{on $M$},
$$
for some positive constant $Q$. Then $v>0$ on $M$ if it is not
identically zero.
\endproclaim

\demo{Proof} Due to the homogeneity of the Paneitz-Branson
equation, we may simply assume $Q=1$ in the following. According
to the above discussion, we have
$$
(-\Delta)^2 w = w^\frac {n+4}{n-4} \tag 2.21
$$
on $\hat \Omega \subset R^n$, where $w = (v\circ\phi\circ
\psi)\hat \eta$. Let us denote $-\Delta w = u$. If we can show
that $u \geq 0$ on $\hat \Omega$, then we easily see that $v>0$ on
$M$ by the strong maximum principle. We will adopt the idea from
the paper of Wei and Xu [WX] for the following argument.

Assume otherwise that $u (0) = -\Delta w (0) < 0$ ( assume $0\in
\hat \Omega$ or after translation if necessary). Notice that $\hat
\Omega = R^n\setminus \psi^{-1}(L(\Gamma))$ where
$\psi^{-1}(L(\Gamma))$ is a compact subset of $R^n$ since $\psi$
is a stereographic projection with respect to a point inside
$\Omega(\Gamma)$. By the assumptions, $w$ is smooth away from
$\psi^{-1}(L(\Gamma))$. Let
$$
\bar f (r) = \frac 1{|S^{n-1}(1)|}\int_{S^{n-1}} f(r\sigma)
d\sigma \tag 2.22
$$
denote the spherical average of the function $f$. Then
$$
\left\{\aligned -\Delta \bar u & = \overline{w^\frac{n+4}{n-4}} \\
                -\Delta \bar w & = \bar u
                \endaligned\right.
$$
at least for very small $r$ and very large $r$. Let us solve the
above equations as follows:
$$
\aligned r^{n-1} \bar u' & = - \int_0^r
s^{n-1}\overline{w^\frac{n+4}{n-4}}ds \\
& = - \frac 1{|S^{n-1}(1)|} \int_0^r \int_{S^{n-1}(1)}
w^\frac{n+4}{n-4}(r\sigma)s^{n-1}d\sigma ds \\
& = - \frac 1{|S^{n-1}(1)|} \int_{B_r(0)} w^\frac{n+4}{n-4} dx.
\endaligned \tag 2.23
$$
Therefore, by above corollary on the integrability, we conclude
that $\bar u$ is differentiable for all $r$, monotonically
decreasing, and
$$
\bar u (r) = u(0) -\int_0^r s^{-n+1}\int_0^s
t^{n-1}\overline{w^\frac{n+4}{n-4}}(t)dtds \tag 2.24
$$
$$
r^{n-1} \bar w' (r) = - \int_0^r s^{n-1} \bar u (s) ds \tag 2.25
$$
and
$$
\bar w (r) = w(0) - \int_0^r s^{-n+1}\int_0^s t^{n-1}\bar u
(t)dtds. \tag 2.26
$$
Plugging our assumption $u(0) < 0$ into (2.24), we have
$$
\bar w(r) \geq \int_0^r s^{-n+1}\int_0^s t^{n-1}\bar u (t)dtds
\geq \frac {-u(0)}{2n} r^2. \tag 2.27
$$
For convenience we denote $q = \frac{n+4}{n-4}$. Applying Jenson's
inequality
$$
\overline{w}^q \leq \overline{w^q}, \tag 2.28
$$
and (2.27) in (2.24), as in [WX], we have
$$
\aligned \bar u (r) & \leq - (\frac {-u(0)}{2n})^q \int_0^r
s^{-n+1} \int_0^s t^{n-1 + 2q}dtds \\
& = - (\frac {-u(0)}{2n})^q \frac 1{n+2q}\frac 1{2q+2} r^{2q+2}
\endaligned
\tag 2.29
$$
and
$$
\aligned \bar w (r) & \geq (\frac {-u(0)}{2n})^q \frac
1{n+2q}\frac 1{2q+2} \int_0^r s^{-n+1}\int_0^t t^{n-1}t^{2q+2}dt
\\
& = (\frac {-u(0)}{2n})^q \frac 1{n+2q}\frac 1{2q+2} \frac
1{n+2+2q}\frac 1{2q+4} r^{2q+4} \\ & \geq (\frac {-u(0)}{2n})^q
\frac 1{(n+4+2q)^4} r^{2q+4}.
\endaligned\tag 2.30
$$
Now the idea is to iterate the above steps over and over. Let
$$
\sigma_k = q\sigma_{k-1} + 4, \sigma_0 = 2.
$$
Then
$$
\aligned \sigma_k & = q^2\sigma_{k-2} + 4q + 4 = q^3\sigma_{k-3} +
4q^2 + 4q + 4 \\ & = q^k \sigma_0 + 4 \frac {q^k- 1}{q-1} \\ & =
2q^k + \frac 4{q-1}q^k - \frac 4{q-1}.\endaligned \tag 2.31
$$
We may calculate and obtain the following after $k$ iterations,
$$
\bar w(r) \geq (\frac {-u(0)}{2n})^{q^k} \prod_{j=1}^k \frac
1{(n+\sigma_j)^{4q^{k-j}}} r^{\sigma_k}. \tag 2.32
$$
To finish the proof we calculate
$$
\prod_{j=1}^k \frac 1{(2\sigma_j)^{4q^{k-j}}} \geq \frac 1 {(4(1+
\frac 2{q-1}))^{\sum_{j=1}^k 4q^{k-j}}} \frac
1{q^{\sum_{j=1}^k4jq^{k-j}}},
$$
where
$$
\sum_{j=1}^k 4q^{k-j} = 4 \frac {q^k-1}{q-1}
$$
and
$$
\sum_{j=1}^k 4jq^{k-j} = 4q^k \sum_{j=1}^k jq^{-j} \leq 4q^k
\int_1^\infty xq^{-x}dx \leq 8q^{k-1}.
$$
Thus
$$
\bar w(r) \geq c_1 (c_2r)^{\sigma_k} \tag 2.33
$$
and $\bar w(r)$ goes to $\infty$ when $k$ goes to $\infty$ for $r
> c_2$, where $c_1$ and $c_2$ are numbers independent of $k$,
which creates contradiction.
\enddemo

\vskip 0.1in\noindent{\bf 3. Convexity}\vskip 0.1in

In this section we use the moving plane method to prove a
convexity result that is similar to what Schoen obtained in [Sc2].
Suppose that $(M, g)$ is a closed locally conformally flat
manifold with positive Yamabe constant but not conformally
equivalent to $(S^n, g_1)$. Suppose that the Kleinian group
$\Gamma$ is the holonomy representation of the fundamental group
of $M$. Let $\Omega(\Gamma)$ be the set of ordinary points. Then
we may assume $(\Omega(\Gamma), \tilde g)$ is a Riemannian
universal cover for $(M, g)$. Let us further assume that the
$Q$-curvature of $(M, g)$ is constant $1$. Then, if let $\tilde g
= u^\frac 4{n-4} g_1$, where $g_1$ the standard round metric on
$S^n$, we have $u > 0$ and
$$
P[g_1] u = u^\frac {n+4}{n-4} \quad\text{in $\Omega(\Gamma)\subset
S^n $}. \tag 3.1
$$
And in this case, we easily see that
$$
u(x) \rightarrow +\infty, \quad\text{as $x \rightarrow L(\Gamma) =
S^n\setminus \Omega(\Gamma)$}. \tag 3.2
$$
Because, as a Riemannian universal cover for $(M, g)$,
$(\Omega(\Gamma), u^\frac 4{n-4} g_1)$ is complete. Using
stereographic projection with respect to some point $p_0 \in
\Omega(\Gamma)$
$$
\psi: \hat \Omega \subset R^n \rightarrow \Omega(\Gamma)\subset
S^n
$$
we may work in the coordinates so that
$$
(\hat \Omega, v^\frac 4{n-4} g_0) = (\Omega(\Gamma)\setminus
\{p_0\}, u^\frac 4{n-4} g_1),
$$
where
$$
v = (u\circ \psi (x))(\frac 2{1+|x|^2})^\frac {n-4}2 > 0. \tag 3.3
$$
Hence
$$
(-\Delta)^2 v = v^\frac {n+4}{n-4} \quad \text{in $\hat
\Omega\subset R^n$}, \tag 3.4
$$
which is equivalent to the system
$$
\left\{\aligned -\Delta w & = v^\frac {n+4}{n-4} \\
                -\Delta v & = w\endaligned\right. \quad\text{in
$\hat \Omega \subset R^n $}, \tag 3.5
$$
We observe

\proclaim{Lemma 3.1} Suppose $(\Omega, u^\frac 4{n-4}g_1)$ is the
universal cover of $(M, g)$. Then
$$
v(x) \rightarrow +\infty, \quad\text{as $x \rightarrow
\psi^{-1}(L(\Gamma)) = R^n\setminus \hat\Omega$}. \tag 3.6
$$
In addition suppose that the scalar curvature for $(M, g)$ is
positive. Then
$$
-\Delta v (x) \rightarrow +\infty, \quad\text{as $x \rightarrow
\psi ^{-1}( L(\Gamma))= R^n\setminus \hat \Omega$}. \tag 3.7
$$
\endproclaim

\demo{Proof} (3.6) is a direct consequence of (3.2). To prove
(3.7) we simply notice that
$$
-\Delta v^\frac {n-2}{n-4} = \hat R v^\frac {n+2}{n-4}
\quad\text{in $\hat\Omega\subset R^n$} \tag 3.8
$$
and
$$
-\Delta v^\frac {n-2}{n-4} = \frac {n-2}{n-4} v^\frac 2{n-4} (
-\Delta v) - \frac {n-2}{n-4}\frac 2{n-4} v^{-\frac {n-6}{n-4}}
|\nabla v|^2, \tag 3.9
$$
where $\hat R = R\circ \phi\circ \psi (x)$. Therefore
$$
-\Delta v = \frac {n-4}{n-2} \hat R v^\frac n{n-4} + \frac 2{n-4}
v^{-1}|\nabla v|^2 \geq \frac {n-4}{n-2} R v^\frac n{n-4}. \tag
3.10
$$
Thus (3.7) holds.
\enddemo

Now we are ready to state and prove a convexity result.

\proclaim{Theorem 3.2} Suppose that $(M^n, g)$ is a closed locally
conformally flat manifold of dimension higher than $4$ which is
not conformally equivalent to $(S^n, g_1)$. And suppose that the
scalar curvature of $(M, g)$ is positive and the $Q$-curvature of
$(M, g)$ is constant $1$. Let $(\Omega(\Gamma), \tilde g)$ be the
Riemannian universal cover for $(M, g)$ associated with the
developing map. Then any round ball $B\subset
\Omega(\Gamma)\subset S^n$ is geodesically convex with respect to
the metric $\tilde g$.
\endproclaim

\demo{Proof} Let $B\subset \Omega(\Gamma)$ be any round ball on
$S^n$. Let $q \in \partial B$ be any given point on the sphere.
Notice that the surface $\partial B$ is umbilical with respect to
the round metric as well as to the metric $\tilde g$ since $\tilde
g$ is conformal to $g_1$. We use stereographic projection $\psi$
with respect to $q$ to construct a coordinate system in which
$\partial B$ is the hyperplane at the origin, {\it i.e.}
$\{x=(x_1, x_2, \cdot, x_n)\in R^n: x_n=0 \}$, and the $x_n$-axis
is in the direction of inward normal direction. Notice that we
have that $\psi^{-1}(L(\Gamma))$ is located $R^n$ with $x_n <0$.
Hence to prove convexity of $B$ with respect to the metric $\tilde
g = v^\frac 4{n-4}g_0$ is equivalent to showing that
$$
\frac {\partial v}{\partial x_n} |_{x_n=0} < 0. \tag 3.11
$$
We will adopt the approach from [Sc2]. Instead of the scalar
curvature equation in [Sc2], we consider the system of elliptic
partial differential equations (3.5). First recall from (3.3) that
$$
v(x) = (u\circ \psi (x)) (\frac 2{1+|x|^2})^\frac {n-4}2,
$$
where $u$ is a smooth function on $\Omega(\Gamma)\subset S^n$ and
therefore is smooth at $\psi (\infty) = q \in \partial B \subset
\Omega(\Gamma)$, thus we have the following expansion
$$
v (x) = \frac 1{|x|^{n-4}}(a_0 + a_i \frac {x_i}{|x|^2} +
a_{ij}\frac{x_ix_j}{|x|^4} +o(|x|^{-2}),  \tag 3.12
$$
where $a_0 > 0$. Hence
$$
\frac{\partial v}{\partial x_n} = - \frac 1{|x|^{n-2}} ((n-4) a_0
x_n - a_n + (n-2) \frac {a_ix_i}{|x|^2}x_n + O(|x|^{-1}), \tag
3.13
$$
$$
w (x) = -\Delta v(x) = 2(n+2)a_0 |x|^{-n+2} + O(|x|^{-n}) \tag
3.14
$$
and
$$
\frac{\partial w}{\partial x_n} = - 2(n+2)(n-2) a_0 x_n |x|^{-n} +
O(|x|^{-n-1}). \tag 3.15
$$
Therefore there exist $C_0$ and $C_1$ such that
$$
\left\{\aligned \frac {\partial v}{\partial x_n} & < 0 \\
                \frac {\partial w}{\partial x_n}  & < 0 \endaligned
                \right. \quad \text{on the set} \ \{x\in
R^n: x_n \geq C_0|x|^{-1} \& \ |x|\geq C_1\}. \tag 3.16
$$
In the following we will use the moving plane method to conclude
(3.11) from (3.16). For $\Lambda \in R$, let
$$
\Sigma_\Lambda = \{x\in R^n: x_n > \Lambda\} \ \text{and} \
S_\Lambda = \{x\in R^n: x_n = \Lambda\}. \tag 3.17
$$
We consider the reflection with respect to the hyperplane $x_n =
\Lambda$
$$
x^\lambda = (x_1, x_2, \cdots 2\Lambda - x_n) \tag 3.18
$$
and define
$$
\left\{\aligned  v^\Lambda (x) & = v(x^\Lambda) \\ w^\Lambda (x) &
= w (x^\Lambda).\endaligned\right.\tag 3.19
$$
Then it is easily seen that
$$
\left\{\aligned -\Delta w^\Lambda (x) & = (v^\Lambda(x))^\frac
{n+4}{n-4} \\ -\Delta v^\Lambda (x) & = w^\Lambda (x)
\endaligned\right. \quad \text{on $\hat\Omega^\Lambda$}. \tag 3.20
$$
Similar to the proofs of Lemma 4.1 and Lemma 4.2 in [GNN] one
knows that there exists $\Lambda_0$ such that, for each $\Lambda >
\Lambda_0$,
$$
\left\{\aligned v^\Lambda (x) & > v(x)\\
                w^\Lambda (x) & > w(x) \endaligned\right. \forall
                x \in \Sigma_\Lambda \setminus
                (\psi^{-1}(L(\Gamma)))^\Lambda. \tag 3.21
$$
Because of Lemma 3.1. Now let
$$
\Lambda^* = \inf\{\Lambda: \ \text{such that (3.21) holds}\}. \tag
3.22
$$
Clearly $\Lambda^*\leq \Lambda_0$. Notice that, applying Hopf
lemma for the elliptic system (cf. [Q])
$$
\left\{\aligned -\Delta (w^{\Lambda} - w) & = (v^{\Lambda})^\frac
{n+4}{n-4} - v^\frac {n+4}{n-4} \\ -\Delta (v^{\Lambda} - v) & =
w^{\Lambda} - w \endaligned\right. \quad \text{in
$\Sigma_{\Lambda} \setminus (\psi^{-1}(L(\Gamma)))^{\Lambda}.$}
\tag 3.23
$$
gives us
$$
0 < \frac{\partial (v^\Lambda - v)}{\partial x_n} |_{x_n =
\Lambda} = -2 \frac{\partial v}{\partial x_n}|_{x_n=\Lambda},
\forall \Lambda \geq \Lambda^*, \tag 3.24
$$
unless $v^\Lambda \equiv v$ which is possible only when the limit
set $L(\Gamma)$ is empty, i.e. $(M, g)$ is conformally equivalent
to $(S^n, g_1)$. Now, again similar to the proof of Lemma 4.4 in
[GNN], we conclude that $\Lambda^*$ must be negative and the
hyperplane $S_{\Lambda^*}$ should not be stopped before hitting
the singular set $\psi^{-1}(L(\Gamma))$. Thus $\Lambda^* \leq 0$
and (3.11) holds.
\enddemo

\vskip 0.1in\noindent{\bf 4. $L^\infty$ estimates and the proof of
main theorem}\vskip 0.1in

In this section we first use the blow-up method to prove the
$L^\infty$ estimate for the conformal factor $u$, when $g_u =
u^\frac 4{n-4}g$ is a metric of positive scalar curvature and
constant $Q$-curvature $1$. Then we use the positivity result
Theorem 2.3 in the previous section to prove the $L^\infty$
estimate for $1/u$. Higher order estimates follow from standard
estimates for linear elliptic partial differential equations.

\proclaim{Theorem 4.1} Suppose that $(M, g)$ is a closed locally
conformally flat manifold of dimension higher than $4$ which is
not conformally equivalent to $(S^n, g_1)$. Then there exists some
constant $C>0$ such that, if $g_u = u^\frac 4{n-4}g$ is a metric
of positive scalar curvature and constant $Q$-curvature $1$, then
$$
\|u\|_{L^\infty(M)} \leq C. \tag 4.1
$$
\endproclaim

\demo{Proof} Suppose otherwise we have a sequence of $u_i$ and a
sequence of points $p_i$ on $M$ such that $(M, u_i^\frac 4{n-4}
g)$ has positive scalar curvature and $Q$-curvature equal to $1$
and
$$
u_i(p_i) = \max_{p\in M} u_i(p) \rightarrow +\infty. \tag 4.2
$$
Let $s_i\in F\subset \Omega(\Gamma)\subset S^n$ such that $\phi
(s_i) = p_i$, where $F$ is a fixed fundamental domain in
$\Omega(\Gamma)$ for the Kleinian group $\Gamma$ and $\phi:
\Omega(\Gamma) \rightarrow M$ is the covering map associated with
the developing map. Without loss of generality, we may assume that
$s_i \rightarrow s \in \bar F\subset \Omega(\Gamma)$. Then the
universal covers $(\Omega(\Gamma), \tilde u_i^\frac 4{n-4} g_1)$
for $(M, u_i^\frac 4{n-4} g)$ satisfy the assumptions in Theorem
3.2 in the previous section. We choose a stereographic projection
$$
\psi: \hat \Omega \subset R^n\rightarrow \Omega(\Gamma)\subset S^n
$$
with respect to some appropriate point inside $\Omega(\Gamma)$ and
work in the coordinates of $(\hat \Omega, \hat u_i^\frac 4{n-4}
g_0)$ for $(\Omega(\Gamma), \tilde u_i^\frac 4{n-4} g_1)$. It is
easily seen that there exists a $\delta >0$ such that $B_\delta
(x_i) \subset \hat \Omega$, where $\phi(\psi(x_i)) = p_i$. Let
$$
v_i (x) = \frac 1{u_i(p_i)} \hat u_i ( x_i + \frac x{u_i^\frac
2{n-4} (p_i)}), \forall x \in B_{u_i^\frac 2{n-4} (p_i)\delta}(0).
\tag 4.3
$$
One may verify that
$$
(-\Delta)^2 v_i = v_i^\frac {n+4}{n-4}, \quad \text{on} \
B_{u_i^\frac 2{n-4} (p_i)\delta}(0) \tag 4.4
$$
and $0 < v_i\leq 1$. Then standard elliptic estimates allow us to
pick up a subsequence $v_{i_k}$ which converges to $v$ in the
$C^{4, \alpha}_{\text{loc}}(R^n)$ topology, where
$$
(-\Delta)^2 v = v^\frac {n+4}{n-4} \quad \text{on} \ R^n \tag 4.5
$$
and $v\geq 0$. Because $-\Delta v \geq 0$ on $R^n$, which was
proved in the proof of Theorem 3.1 in [WX], and $v (0) = 1$, $v$
is actually positive. By a classification theorems in [WX], we
know that $(R^n, v^\frac 4{n-4}g_0)$ is the round sphere. Hence, a
Euclidean ball $B_K(0)$ of sufficiently large radius is strictly
geodesically concave in the round metric $v^\frac 4{n-4}g_0$.
Therefore the balls $B_K(0)$ are geodesically concave in the
metrics $v_{i_k}^\frac 4{n-4} g_0$ when $k$ is sufficiently large.
That is to say the Euclidean ball
$B_{u_{i_k}^{-1}(p_{i_k})K}(x_{i_k})$ is geodesically concave in
the metric $\hat u_i^\frac 4{n-4} g_0$, which is a contradiction
with Theorem 3.2 in the previous section. Thus the proof is
completed.
\enddemo

Now we state and prove our main theorem in this note.

\proclaim{Theorem 4.2} Let $(M, g)$ be a closed locally
conformally flat manifold of dimension greater than $4$ with
positive Yamabe constant and Poincar\'{e} exponent less than
$\frac {n-4}2$. Suppose that $(M, g)$ is not conformally
equivalent to $(S^n, g_1)$ and that the Paneitz constant $P(M,
[g])$ is positive. Then there exists $C_k
>0$ for each $k\in N$ such that, if $g_u = u^\frac 4{n-4}g$ is of
positive scalar curvature and constant $Q$-curvature $1$, then
$$
\|u\|_{C^k(M)} + \|\frac 1u \|_{C^k(M)} \leq C_k. \tag 4.6
$$
\endproclaim
\demo{Proof} By the previous theorem and standard elliptic theory
we only need to show that there is a uniform positive lower bound
for the function $u$, if $g_u = u^\frac 4{n-4}g$ is of positive
scalar curvature and constant $Q$-curvature $1$. Suppose
otherwise, there is a sequence of $u_i$ such that
$$
\min_{p\in M} u_i(p) \rightarrow 0. \tag 4.7
$$
Then, by the previous theorem and standard elliptic theory, we may
take a subsequence $u_{i_k}$, which converges to $u \geq 0 $ in
$C^{4, \alpha}(M)$, where
$$
P[g] u = u^\frac {n+4}{n-4} \quad \text{on} \ M
$$
and $u (p) =0 $ for some $p\in M$. But the existence of such $u$
contradicts with our Theorem 2.3 in Section 2, unless $u\equiv 0$.
If indeed $u\equiv 0$, then we would have
$$
\frac {\int_M (Qdv)[g_{u_{i_k}}]}{(\int_M dv[g_{u_{i_k}}])^\frac
{n-4}n} = (\int_M u_{i_k}^\frac {2n}{n-4}dv_g)^\frac 4n
\rightarrow 0 \quad\text{as $k\rightarrow\infty$}, \tag 4.8
$$
which contradicts with the assumption that the Paneitz constant
$P(M, [g])$ is positive. So the proof is finished.
\enddemo

\vskip 0.1in \noindent {\bf References}:

\roster \vskip0.1in
\item"{[Au]}" T. Aubin, The scalar curvature, differential geometry
and relativity, (Cahen and Flato, eds.), Reidel, Dordrecht 1976.

\vskip0.1in
\item"{[Br]}" Thomas Branson, Group representations arising from
Lorentz conformal geometry, J. Func. Anal. 74 (1987), 199-291.



\vskip 0.1in
\item"{[CHY]}" S.-Y. A. Chang, F. Hang and P. Yang, On a class of
locally conformally flat manifold, Preprint 2003.


\vskip0.1in
\item"{[CQY2]}" S.-Y. A. Chang, J. Qing and P. Yang, On the finiteness
of Kleinian groups in general dimensions, J. Reine Angew. Math.,
571 (2004), 1-17.

\vskip0.1in
\item"{[CY]}" S.-Y. A. Chang and P. Yang, Nonlinear differential
equations in conformal geometry, Proceedings of ICM 2002.

\vskip0.1in
\item"{[DHL]}" Z. Djadli, E. Hebey and M. Ledoux, Paneitz type
operators and applications, Duke Math. J. 104(2000) no. 1,
129-169.

\vskip0.1in
\item"{[DMA]}" Z. Djadli, A. Malchiodi and M. Ahmedou, Prescribing
a fourth order conformal invariant on the standard sphere, Part II
- blow-up analysis and applications, Preprint 2001.

\vskip0.1in
\item"{[GNN]}" B. Gidas, W.-M. Ni and L. Nirenberg, Symmetry and
related properties via maximum principle, Comm. Math. Phys. 68
(1979), no. 3, 209-243.

\vskip0.1in
\item"{[G]}" Maria Gonzalez, Singular sets of a class of locally
conformally flat manifolds, Preprint 2004.

\vskip 0.1in
\item"{[HR]}" E. Hebey and E. Robert, Compactness and global estimates for
the geometric equation in higher dimensions, E.R.A./A.M.S.,
Vol. 10, pp. 135-141, 2004.


\vskip0.1in
\item"{[P]}" S. Paneitz, A quadratic conformally covariant
differential operator for arbitrary pseudo-Riemannian manifolds,
Preprint 1983.

\vskip0.1in
\item"{[Q]}" Jie Qing, A Priori Estimates for Positive
Solutions of Semi-linear Elliptic Systems. J. Partial Differential
Equations, Vol.1, No.2, Series A(1988), pp 61-70.

\vskip0.1in
\item"{[Sc1]}" Rick Schoen, Conformal deformation of a Riemmanian
metric to constant scalar curvature, J. Diff. Geom. 6 (1984),
479-495.

\vskip0.1in
\item"{[Sc2]}" Rick Schoen, On the number of constant scalar
curvature metrics in a conformal class, Differential Geometry,
311-32o, Pitman Monogr. Surveys Pure Appl. Math., 52, Longman Sci.
Tech., Harlow, 1991.

\vskip0.1in
\item"{[SY]}" R. Schoen and S.T. Yau, Conformally flat manifolds,
Kleinian groups, and scalar curvature, Invent. Math. 92 (1988),
47-71.

\vskip0.1in
\item"{[Tr]}" N. Trudinger, Remarks concerning the conformal
deformation of Riemannian structures on compact manifolds, Ann.
Scuola Norm. Sup. Pisa Cl. Sci (3) 22 (1968), 265-274.

\vskip0.1in
\item"{[WX]}" Juncheng Wei and Xingwang Xu, Classification of
solutions of higher order conformally invariant equations, Math.
Ann. 313 (1999), 207-228.

\vskip0.1in
\item"{[Ya]}" H. Yamabe, On a deformation of Riemannian structures
on compact manifolds, Osaka J. Math. 12 (1960), 21-37.

\endroster

\enddocument